\documentclass{amsart}
\usepackage[shortalphabetic]{amsrefs}
\usepackage[utf8]{inputenc}
\usepackage{amsthm}
\usepackage[utf8]{inputenc}
\usepackage{appendix}
\usepackage{cite}
\usepackage{amssymb,amsmath,amsthm}
\usepackage{hyperref}
\newcommand{\amsprimary}[1]{{\footnotesize\noindent AMS 2010 \textit{Mathematics subject
classification:} Primary #1\vspace{1pc}}}
\newcommand{\keywordsnames}[1]{{\footnotesize\noindent\textit{Key words:} #1\vspace{1pc}}}

\newtheorem{theorem}{Theorem}
\newtheorem{teo}{Theorem}

\newtheorem{corollary}[teo]{Corollary}

\theoremstyle{definition}

\theoremstyle{remark}

\title[Dirichlet problem in cones]{An Observation on the Dirichlet problem at infinity in Riemannian cones}
\author{Jean C. Cortissoz}
\email{jcortiss@uniandes.edu.co}
\address{Department of Mathematics, Universidad de los Andes, Bogot\'a DC, Colombia}
\date{}

\begin{document}

\maketitle

\begin{abstract}
    In this short paper we show a sufficient condition for the solvability of
    the Dirichlet problem at infinity in Riemannian cones (as defined below).
    This condition is related to a celebrated result of Milnor that
    classifies parabolic surfaces. When applied to 
    smooth Riemannian manifolds with a special type of metrics (which generalise
    rotational symmetry) we obtain generalisations of classical criteria for the
    solvability of the Dirichlet problem at infinity. Our proof is short and elementary:
    it uses separation of variables
    and comparison arguments for ODE's.
\end{abstract}

\keywordsnames{Dirichlet problem, bounded harmonic function, Riemannian cone.}

{\amsprimary {31C05, 53C21}}

\tableofcontents

\section{Introduction}

We define a Riemannian cone as follows. Let $N$ be an $\left(n-1\right)$-dimensional compact manifold, and consider the
quotient space
\[
M=\left(N\times\left[0,\infty\right)\right)/\left(N\times\left\{0\right\}\right),
\]
endowed with a complete metric that can be written as
\begin{equation}
\label{eq:warped_metric}
g=dr^2+\phi\left(r\right)^2 g_{N},
\end{equation}
where $g_{N}$ is any smooth metric on $N$. 
We shall assume that $\phi$ is smooth, and that $\phi\left(0\right)=0$ and $\phi'\left(0\right)=1$.
It is usual to define the cone metric using $g=r^2$, so our definition of a Riemannian cone is a bit
more general. Also, if $N=\mathbb{S}^{n-1}$ and $g$ is smooth up to the pole (or vertex), with
$g_{N}$ any smooth metric on the unit $\left(n-1\right)$-dimensional sphere, this definition
of a cone defines a family of metrics in $\mathbb{R}^n$ that includes those with rotational symmetry,
which are those obtained when $N=\mathbb{S}^{n-1}$ and it is endowed with the round metric.

\medskip
The Laplacian in a cone can be written for $r>0$ as:
\[
\Delta=\frac{\partial^2}{\partial r^2}+\left(n-1\right)\frac{\phi'}{\phi}\frac{\partial}{\partial r}
+\frac{1}{\phi^2}\Delta_N,
\]
where $\Delta_N$ is the Laplacian in $N$, and as expected, a function
$u:M\longrightarrow \mathbb{R}$ if $\Delta u=0$.
We will prove the following result -for applications see below.
\begin{theorem}
\label{thm:main}
Assume that
\begin{equation}
\label{ineq:fundamental}
\int_1^\infty \frac{1}{\phi\left(s\right)}\,ds<\infty,
\end{equation}
and, if $\mbox{dim}\left(N\right)\geq 2$, that there is an $R_0$ such that if $r>R_0$ then 
$\phi'\left(r\right)\geq 0$.
Then for any $f:N\longrightarrow \mathbb{R}$ regular enough
(to avoid technicalities just set $f\in C^{\infty}\left(N\right)$), the Dirichlet problem
is uniformly solvable at infinity. By this we mean that there is
a harmonic function $u:M\longrightarrow \mathbb{R}$ such that
\[
\lim_{r\rightarrow\infty}u\left(r,\omega\right)=f\left(\omega\right)
\quad\mbox{uniformly}.
\]
\end{theorem}

The definition of uniform solvability at infinity we used in the previous statement is
stronger than the one commonly used when studying the 
Dirichlet problem at infinity in Cartan-Hadamard manifolds, i.e., when both definitions
apply, uniform solvability implies solvability in the sense of Choi \cite{Choi}.
By smooth enough we mean that $f$ must be regular enough so that its 
expansion in eigenfunctions of the Laplacian $\Delta_N$
converges uniformly. The fact that enough regularity of $f$ implies
the absolute convergence of its expansion in eigenfunctions of the Laplacian
 was shown by Peetre in \cite{Peetre}.
Theorem \ref{thm:main} implies that cones with metrics satisfying (\ref{ineq:fundamental}) 
have a wealth
of nontrivial bounded harmonic functions. 

\medskip
For applications of our main result, we specialize to the 
case when $N=\mathbb{S}^{n-1}$ with an arbitrary smooth metric so that
the cone is a smooth Riemannian manifold.
From Theorem \ref{thm:main}, given $f\in C\left(N\right)$, taking a sequence
$f_n$ of smooth functions such that $f_n\rightarrow f$ uniformly, then solving
the Dirichlet problem for each $f_n$ (in this case the solution is unique
by the asymptotic maximum principle of Choi, Proposition 2.5 in \cite{Choi}), 
we obtain a uniformly convergent sequence of harmonic functions 
(again by the asymptotic maximum principle), and
thus:

\begin{corollary}
\label{cor:main}
Assume that (\ref{ineq:fundamental}) holds and, if $\mbox{dim}\left(N\right)\geq 2$, that there is an $R_0\geq 0$ such that if $r>R_0$ then $\phi'\left(r\right)\geq 0$.
Then for any continuous $f:N\longrightarrow \mathbb{R}$ the Dirichlet problem
is uniquely solvable at infinity. By this we mean that there is
a harmonic function $u:M\longrightarrow \mathbb{R}$ such that
\[
\lim_{r\rightarrow\infty}u\left(r,\omega\right)=f\left(\omega\right)
\quad\mbox{in the cone topology (see Section \ref{section:preliminaries})}.
\]
\end{corollary}

\medskip
Corollary \ref{cor:main} seems to be new: notice that
we do not require the manifold to be of nonpositive curvature, that is, to 
be Cartan-Hadamard; however, if we assume
the manifold to be Cartan-Hadamard, the hypothesis $\phi'$ can be dropped,
as it is automatically satisfied. 
From the previous corollary we obtain the following.
\begin{corollary}
\label{cor:main2}
Let $g$ be a metric of the form (\ref{eq:warped_metric}). If there is
an $\epsilon>0$ such that $-\phi''/\phi\leq -\left(1+\epsilon\right)/r^2\log r$,
for large enough $r$, and $\phi$ is unbounded,
then for any $f\in C\left(N\right)$, the Dirichlet problem is uniquely solvable at infinity.
\end{corollary}

The proof of this corollary is as follows. By Milnor's argument in \cite{Mil}, we have that
\[
\int_1^{\infty}\frac{1}{\phi\left(s\right)}\,ds<\infty.
\]
The condition $\phi'\geq 0$ is automatically satisfied
for large enough $r$: if $\phi'<0$ for all $r$ large enough, then
$1/\phi$ would be increasing for $r$ large enough and then it would not
be integrable, and hence at some point, say $r_1$, we must 
have $\phi'\left(r_1\right)\geq 0$, and since the radial curvature $-\phi''/\phi$ is negative,
$\phi'$ will be nonnegative from then on; thus the hypotheses of Corollary \ref{cor:main} hold
and Corollary \ref{cor:main2} follows.
Again, notice that
Corollary \ref{cor:main2} includes metrics that are {\bf not necessarily rotationally symmetric},
since in the case of $N=\mathbb{S}^{n-1}$ it is not required 
that the metric it carries to
be the round metric nor conformal to it, in contrast with the results proved
in \cite{Choi} and \cite{March} which require either one 
(specially \cite{Choi} where rotational symmetry is strongly required
in the calculations -see Section 3. On the other hand, in \cite{March} 
rotational symmetry does not seem to be essential, in the sense that 
the metric that $\mathbb{S}^{n-1}$ carries does not need to be the round metric,
but the author only proves the transience of the manifold). In this 
sense, Corollary \ref{cor:main2} is new for dimensions greater or equal than 4.
As an aside comment: when Milnor's criteria is used, the "$\phi$ unbounded" part
in its statement is usually replaced by saying that
the manifold has everywhere nonpositive curvature, and so, again, the assumption
of $\phi$ being unbounded might be replaced by assuming the stronger condition that the manifold
is a Cartan-Hadamard manifold, which is the assumption made by Choi in \cite{Choi}.

\medskip
Regarding (\ref{ineq:fundamental}), in a beautiful work \cite{Neel}, R. Neel shows that
in Cartan-Hadamard surfaces if we write the metric as
\begin{equation}
\label{metric:general}
    dr^2+J\left(r,\theta\right)^2d\theta^2,
\end{equation}
and 
\begin{equation}
\label{ineq:transience}
\int_1^{\infty}\frac{1}{J\left(r,\theta\right)}\,dr<\infty,
\end{equation}
then the Dirichlet problem at infinity is uniquely solvable. His approach
is probabilistic,  and his result
strengthens the following result of Doyle, at least
in the case of surfaces, which is proved in \cite{Grigoryan}
: For a Riemannian manifold with a metric written as (\ref{metric:general}),
 (\ref{ineq:transience}) implies transience.

\medskip
Also, our proof is quite elementary, perhaps strikingly simple 
when compared with the published proofs of the results mentioned above, 
and it reveals that there is a notion of weak solvability
for the Dirichlet problem at infinity
that has not been treated before in the literature (to the best of our knowledge),
and which perhaps deserves more consideration. In fact, we have the following:
\begin{theorem}
Assume that there is an $R_0$ such that if $r>R_0$ $\phi'\left(r\right)\geq 0$ and that
it satisfies (\ref{ineq:fundamental}).
Then for any $f\in L^{2}\left(N\right)$ 
the Dirichlet problem
weakly is solvable at infinity. By this we mean that there is function
$u:M\longrightarrow \mathbb{R}$ which satisfies $\Delta u=0$ in the
sense of distributions (and hence classically when $M$ is a manifold), and such that
\[
\lim_{r\rightarrow}u\left(r,\omega\right)=f\left(\omega\right) \quad \mbox{in}
\quad L^2\left(N\right).
\]
\end{theorem}

\medskip
The main ingredient in our proof of Theorem \ref{thm:main} is separation of variables
(for a precursor of our proof see \cite{Dodziuk}). We shall show that we can write a 
whole family of harmonic functions as
\[
u\left(r,\omega\right)=\sum_m \varphi_{m}\left(r\right)\left(\sum_k c_{m,k} f_{m,k}\left(\omega\right)\right),
\]
where $f_{m,k}$ are the eigenfunctions of the Laplacian $\Delta_{g_{N}}$ of $N$
with eigenvalue $\lambda_m^2$. 
Then we show that that the $\varphi_m$ can be chosen to be nonnegative, increasing and bounded.
In the case when the sectional curvature satisfies 
$K\leq -1$ everywhere, we can even have for large enough $r$
\begin{equation}
\label{ineq:eigenfunctions}
0\leq \varphi_m\left(r\right)\leq A_m\tanh^{\lambda_m} r,
\end{equation}
for a convenient constant $A_m$.

\medskip
Since the eigenfunctions of the Laplacian in $N$ form a basis for
$L^2\left(N\right)$, then we can solve the Dirichlet problem.

\medskip
The Dirichlet problem at infinity has
a rich history full of deep and interesting results (see for instance \cite{A, AS, Choi, Ji, Sullivan}).
Theorem \ref{thm:main} and its consequences represent an improvement in the case of rotationally symmetric metrics 
(by also including a bit more general type of metrics)
as given in \cite{Choi}, and it is a natural extension of the classical result of Milnor in 
\cite{Mil}.
Also, our main result bears some resemblance with that of March in \cite{March}, where under the
hypothesis
\begin{equation}
\label{ineq:march}
\int_1^{\infty} \phi\left(r\right)^{n-3}\int_r^{\infty}\phi\left(\rho\right)^{1-n}\,d\rho\,dr<\infty,
\end{equation}
and rotational symmetry, 
the author proves the existence of nonconstant bounded harmonic functions in $M$. 
The reader will find
the proof of Theorem \ref{thm:main}, our main result, in Section \ref{section:proof_main}.

\section{Preliminaries}
\label{section:preliminaries}

Here we define what is to be understood as to solve the Dirichlet problem at infinity. A way
of doing it is by firstly defining a
a compactification. To this end, define the 
set $\overline{M}$ 
\[
\overline{M}=N\times \left[0,\infty\right]/\left(N\times\left\{0\right\}\right),
\]
where $\left[0,\infty\right]$ is a compactification of $\left[0,\infty\right)$
which is homeomorphic to $\left[0,1\right]$.
The subspace $\partial_{\infty}M:=N\times\left\{\infty\right\}$ 
plays the role of
a boundary, and in fact, when
$N$ is homeomorphic to $\mathbb{S}^{n-1}$, $\overline{M}$ has the strcuture
of a topological manifold with boundary.
Given the previous definition, a way of defining that {\bf the Dirichlet problem is solvable at
infinity} is as follows: Given $f\in C\left(N\right)$ there exists a function $u:\overline{M}\rightarrow \mathbb{R}$
which is in $C\left(\overline{M}\right)$, is harmonic in $M$, and such that its restriction
to $\partial_{\infty}M$ is $f$. In the case of Cartan-Hadamard manifolds, this definition
of solvability coincides with the definition of solvability given by Choi, who uses the {\bf cone topology}
as defined by Eberlein-O'Neill \cite{Eberlein}, and which is equivalent,
in the sense of homeomorphism, to the one defined for the
compactification above: the resulting spaces in both cases are homeomorphic to the closed $n$-ball.

\medskip
A stronger definition of
solvability was used in the
statement of the Theorem \ref{thm:main}, let us recall it. 
{\bf We shall say that the Dirichlet problem is uniformly solvable at infinity}
if there is a harmonic function $u:M\longrightarrow \mathbb{R}$ such that
\[
\lim_{r\rightarrow \infty}u\left(r,\omega\right)=f\left(\omega\right)
\]
uniformly on $\omega\in N$. Notice that using this definition, if $f\in C\left(N\right)$ then we can extend
$u$ to $\overline{M}$ continuously and hence this definition of solvability 
implies the one given above. 
In the case of $N=\mathbb{S}^{n-1}$ (not necessarily with the round metric) uniform 
solvability implies solvability in the sense of Choi in \cite{Choi}. %It must be pointed out
%that Choi's definition of solvability is more akin to the notion of solvability of the
%Dirichlet problem in bounded domains.

\medskip
Observe then that with the definitions of solvability given 
in the previous paragraph, if the Dirichlet problem is solvable at infinity
for a given continuous data, then the corresponding harmonic
extension is bounded; so solving 
the Dirichlet problem at infinity gives a method for proving the existence of bounded nonconstant harmonic functions.

\medskip

Also, we can define the concept of weak solvability at infinity
for the Dirichlet problem. In this case, given $f\in L^{2}\left(N\right)$
we say that the Dirichlet problem is weakly solvable at infinity with
boundary data $f$ if there is a function $u:M\longrightarrow \mathbb{R}$
which is harmonic in the sense of distributions (and thus harmonic by Weyl's lemma), and such that
\[
\lim_{r\rightarrow \infty} u\left(r,\omega\right)= f\left(\omega\right)
\quad 
\mbox{in} \quad L^{2}\left(N\right).\]

\section{A proof of Theorem \ref{thm:main}}
\label{section:proof_main}

The starting point of our proof of Theorem \ref{thm:main}
is a simple computation. We use separation of variables to find the 
equation the functions $\varphi_m$ must satisfy so that the product
$\varphi_m\left(r\right) f_{m,k}\left(\Omega\right)$ is harmonic.
We let
 $f_{m,k}$ be the eigenfunctions of the $m$-th eigenvalue, $\lambda_m^2$
 (we use the convention that $\lambda_m\geq 0$), of the
Laplacian of $N$.
It is elementary to prove then that
the equation to be satisfied by the $\varphi_m's$ is
\begin{equation}
\label{eq:radialeigenfunctions}
\varphi_m''+\left(n-1\right)\frac{\phi'}{\phi}\varphi_m'-\frac{\lambda_m^2}{\phi^2}\varphi_m=0.
\end{equation}
First we show that there is a solution such that $\varphi_m\left(0\right)=0$, if $m\neq 0$, and that we can choose
$\varphi_m>0$ near 0 (for $m=0$ we just choose the constant function 1). 
Indeed, since $\phi\left(r\right)\sim r$, $r=0$ is a regular singular point
of the equation, and thus, near $r=0$ a solution can be written as 
\[
q\left(r\right)=r^{k}p\left(r\right),
\]
where $k=\dfrac{-\left(n-2\right)+\sqrt{\left(n-2\right)^2+4\lambda_m^2}}{2}>0$ satisfies the indicial equation
\[
k\left(k-1\right)+\left(n-1\right)k-\lambda_m^2=0,
\]
where $p$ is smooth and $p\left(0\right)=1$ (see page 45 in \cite{Bocher}, a classical paper
by M. Bôcher, and Chapters
4 and 5 in \cite{Levinson}). From this our assertion follows.

\medskip
Next we show that $\varphi_m$ is nondecreasing. From its general form,
it is clear that $\varphi_m'\left(0\right)> 0$ near $0$. Assume then that at some point
$\varphi_m'=0$ occurs for the first time. Then, since up to that 
first point $\varphi_m\geq 0$, using equation
(\ref{eq:radialeigenfunctions}), shows that $\varphi_m''\geq 0$, which in turn
implies our claim.

\medskip
Define the function
\[
\eta_m\left(r\right)=\exp\left[-\int_{r}^{\infty}\frac{\lambda_m}{\phi}\,ds\right].
\]
Of course here is where we need the fundamental assumption that
\[
\int_1^{\infty}\frac{1}{\phi}\,ds<\infty.
\]
\medskip
A straightforward computation shows that if $r>R_0$, then
\[
\mathcal{L}_m \eta_m =\left(n-2\right)\frac{\phi'}{\phi}\lambda_m \eta_m \geq 0.
\]
This is why we need the hypothesis $\phi'\left(r\right)\geq 0$
for $r$ large enough.
We note in passing that when $\mbox{dim}\left(N\right)=n-1=1$, the $\eta_m$'s thus defined give
explicit solutions to (\ref{eq:radialeigenfunctions}) (this was pointed out to me
by J. E. Bravo): this is the reason why in dimension 2 the assumption
$\phi'\geq 0$ is not required.

\medskip
Next, we are going to use the inequality above to show that for $R_1>R_0>0$, if 
$r>R_1$, it then holds that $0\leq \varphi_m\left(r\right)\leq A_m\eta_m\left(r\right)$
for an appropriate constant $A_m$.
So let $A_m>0$ large enough so that
\[
\varphi\left(R_1\right)<A_m\eta_m\left(R_1\right),\quad 
\varphi_m'\left(R_1\right)<A_m\eta_m'\left(R_1\right).
\]
We now show that $\left(A_m\eta_m-\varphi_m\right)'\left(r\right)\geq 0$ must
also hold for $r>R_1$. Say the strict inequality holds up to $r_1>R_1$ and that
equality is attained at $r=r_1$. Then, as up to $r=r_1$ we have that $h_m:=A_m\eta_m-\varphi_m>0$,
then we get
\[
h_m''\left(r_1\right)\geq \frac{\lambda_m^2}{\phi^2}h_m\left(r_1\right)> 0,
\]
and hence right after $r_1$ we again have $h_m'>0$. But then, this implies
that $A\eta_m\geq \varphi_m$ right after $r_1$. Therefore, we have that for
$r\geq R_1$ that $\varphi_m\leq A_m \eta_m$, and hence $\varphi_m$ is bounded above.

\medskip
Notice that in the case that the sectional curvature satisfies
$K\leq -1$ everywhere, then, by the Bishop-Gromov theorem
(which in the case of rotational symmetry reduces to a simple ODE 
comparison argument), 
$\phi\left(r\right)\geq \sinh r$, and hence
we have that $0\leq \varphi_m \leq A_m\tanh^{\lambda_m} r$, and claim
(\ref{ineq:eigenfunctions})
holds.
In any case, the previous estimates show that 
by multiplying by appropriate constants, we may assume that
$\lim_{r\rightarrow\infty} \varphi_m\left(r\right)=1$. 

\medskip
Given $f\in C^{\infty}\left(N\right)$ we can represent it as
\[
f\left(\omega\right)=\sum_m \sum_k c_{m,k}f_{m,k}\left(\omega\right).
\]
From this we get a harmonic extension
\[
u\left(r,\omega\right)=\sum_m \varphi_m\left(r\right)\sum_k c_{m,k}f_{m,k}\left(\omega\right).
\]
If $f$ is smooth enough (assume $f\in C^{\infty}\left(N\right)$ from now on), it is not difficult
to show that $u$ is twice differentiable and harmonic. 
All we need to prove is that the harmonic function $u$ satisfies the boundary
conditions at infinity. Let $\epsilon>0$, we estimate as follows. 
By the triangular inequality,
\begin{eqnarray*}
\left|f\left(\omega\right)-u\left(r,\omega\right)\right|&\leq &\sum_m \left(1-\varphi_{m}\left(r\right)\right)
\left|\sum_k c_{m,k} f_{m,k}\left(\omega\right)\right|.
\end{eqnarray*}
Pick $M$ such that 
\[
\sum_{m\geq M} 
\left|\sum_k c_{m,k} f_{k,m}\left(\omega\right)\right|\leq \frac{\epsilon}{2}.
\]
This can be done by Peetre's result as soon as $f$ is smooth enough.
Let $R>0$ be such that if $r\geq R$, for $m=0,1,\dots, M$ 
\[
1-\varphi_m\leq \frac{\epsilon}{2L},
\]
where $L$ bounds 
\[
\sum_{m\leq M} 
\left|\sum_k c_{m,k} f_{k,m}\left(\omega\right)\right|
\]
from above. Under this considerations, we obtain that, for $r\geq R$,
\[
\left|f\left(\omega\right)-u\left(r,\omega\right)\right|< \epsilon,
\]
and our claim is now proved. This finishes the proof of Theorem \ref{thm:main}.

\subsection{Weak Solvability}
Notice that if we only require $f\in L^{2}\left(N\right)$, the harmonic extension
$u$ is harmonic in the sense of distributions (and hence harmonic) and we can also show that
$u\left(r,\omega\right)\rightarrow f\left(\omega\right)$ as $r\rightarrow \infty$ in $L^{2}\left(N\right)$,
and in consequence that the Dirichlet problem is weakly solvable. 
First we show that $u$ is weakly harmonic (and hence almost everywhere equal to a smooth
harmonic function, and thus harmonic). This is standard:
Given $\varphi \in C^{\infty}_0\left(N\right)$ we let
\[
u_l\left(r,\omega\right)=\sum_{m\leq l} \varphi_m\left(r\right)\sum_k c_{m,k}f_{m,k}\left(\omega\right).
\]
As $u_l$ is harmonic, then
\[
\int_M u_l\left(x\right)\Delta\varphi\left(x\right)\,d\mu_g=0.
\]
Since the $\varphi_m$ are uniformly bounded, and $u_l\rightarrow u$ in $L^2_{loc}\left(M\right)$,
then we must have that
\[
\int_M u\left(x\right)\Delta\varphi\left(x\right)\,d\mu_g=0,
\]
which is what we wanted to show.

\medskip
The following computation,
which can be originally found in \cite{Bravo} and that we reproduce for the convenience 
of the reader, shows that the boundary data is satisfied in an $L^2$-sense as described above.
Denoting by $\left\|\cdot\right\|_2$ the $L^2\left(N\right)$-norm,
using the orthogonality of the eigenfunctions,
 we can estimate
\begin{eqnarray*}
\left\|f\left(\omega\right)-u\left(r,\omega\right)\right\|^2_2&=&\left\|\sum_m \left(1-\varphi_{m}\left(r\right)\right)
\sum_k c_{m,k} f_{k,m}\left(\omega\right)\right\|^2_2\\
&=&\sum_m \left|1-\varphi_{m}\left(r\right)\right|^2
\left\|\sum_k c_{m,k} f_{m,k}\left(\omega\right)\right\|^2_2.
\end{eqnarray*}
Pick $M$ such that 
\[
\sum_{m\geq M} 
\left\|\sum_k c_{m,k} f_{m,k}\left(\omega\right)\right\|^2_2\leq \frac{\epsilon^2}{8},
\]
and let $R>0$ such that if $r\geq R$, for $m=0,1,\dots, M$ we have that
\[
\left|1-\varphi_m\right|\leq \frac{\epsilon}{2K},
\]
where $K$ is such that
\[
\sum_{m\leq M} 
\left\|\sum_k c_{m,k} f_{m, k}\left(\omega\right)\right\|^2_2\leq K.
\]
Notice that $\left|1-\phi_m\left(r\right)\right|\leq 1$. Hence, putting all this together, 
we find that given any $\epsilon>0$ 
for $r$ large enough
\[
\left\|f\left(\cdot\right)-u\left(r,\cdot\right)\right\|_2<\epsilon,
\]
which is what we wanted to prove.

\section{Some Remarks}

The case of $N$ being a Lie group, Taylor in \cite{Taylor} gave sufficient
conditions for the eigenfunction expansion of $f$ to converge uniformly to 
$f$. For instance if $N=\mathbb{S}^3$ with the round metric,
then if $f\in C^{\frac{1}{2}}\left(\mathbb{S}^3\right)$, the Dirichlet problem is 
not just solvable but uniformly solvable
at infinity in $\mathbb{R}^4$ endowed with a rotationally symmetric metric such
that the factor $\phi$ satisfies (\ref{ineq:fundamental}) and such that 
$\phi'\geq 0$ (or if $\mathbb{R}^4$ with the given metric is a Cartan-Hadamard manifold). 
On the other hand, if $N$ is a Lie 
group of dimension $4k$, it is only required that $f\in H^{2k}\left(N\right)$ for 
the limit $\lim_{r\rightarrow\infty} u\left(r,\omega\right)=f\left(\omega\right)$ to be uniform.

\subsection{Acknowledgements}
The author wants to thank J. E. Bravo for discussing 
and bouncing some of the ideas used in this paper; 
the staff of \textquotedblleft La Valija de Fuego\textquotedblright  for giving the
author a marvelous space to think about mathematics; my students
for making me think about separation of variables; and specially my wonderful wife and baby son 
for being 
so patient with me when I am \textquotedblleft spaced out\textquotedblright.
%\subsection{Interesting Observations}
%We have thus proved that the weak Dirichlet problem is solvable for the
%metrics considered in this paper. Also, this functions belong to $L^{2}\left(M\right)$.

\end{document}